\documentclass[12pt]{article}
\usepackage{graphicx} 
\usepackage{amsmath,amsthm,amssymb,enumerate,xcolor}

\title{Integrability of the magnetic geodesic flow on the sphere with a constant 2-form}
\author{Alexey V. Bolsinov\footnote{ School of Mathematics,
 Loughborough University,
 LE11 3TU, UK  and Institute of Mathematics and
Mathematical Modeling, Almaty, Kazakhstan\ \ \quad {\tt   A.Bolsinov@lboro.ac.uk}  } \quad \& \quad  Andrey Yu. Konyaev\footnote{Faculty of Mechanics and Mathematics, Moscow State University,  and  Moscow Center for Fundamental and Applied Mathematics, 119992, Moscow Russia, and Institute of Mathematics and
Mathematical Modeling, Almaty, Kazakhstan
 \ \ \quad {\tt  maodzund@yandex.ru}}  \quad \& \quad   Vladimir S.\ Matveev\footnote{
Institut f\"ur Mathematik, Friedrich Schiller Universit\"at Jena,
07737 Jena,  Germany  \ \ \quad {\tt  vladimir.matveev@uni-jena.de}}}

\newtheorem{theorem}{Theorem} 
\newtheorem{remark}[theorem]{Remark}

\newtheorem{fact}[theorem]{Fact}

\newtheorem{corollary}[theorem]{Corollary}
\numberwithin{theorem}{section}

\newcommand{\weg}[1]{}
\newcommand{\dd}{\operatorname{d}}

\usepackage{biblatex} 
\begin{document}

\maketitle
\begin{abstract}  We prove a recent conjecture of Dragovic et al \cite{dragovic}
stating that the magnetic geodesic flow on the standard sphere $S^n\subset \mathbb R^{n+1}$ whose magnetic 2-form is the restriction of a constant 2-form from $\mathbb{R}^{n+1}$
is Liouville integrable. The integrals are quadratic and linear  in momenta.


\end{abstract}

\section{Setup and results}

Recall that by  
{\it magnetic geodesic flow }   on a Riemannian (or pseudo-Riemannian) manifold $(M, g)$ endowed with a  closed differential 2-form $\omega$,  one understands 
the  Hamiltonian system with respect to the ``perturbed'' symplectic form 
\begin{equation}
    \label{eq:pert}
  \Omega_{\mathrm{pert}}:= \omega_{ij}\dd x^i\wedge  \dd x^j +  \sum_{i=1}^n  \dd p_i\wedge \dd x^i 
  \end{equation} generated by the 
Hamiltonian 
\begin{equation}
\label{eq:Hunpert}
H(x,p)= \tfrac{1}{2}\sum_{i,j=1}^n  g^{ij}p_ip_j. 
\end{equation}
Here  $x= (x^1, \dots, x^n)$ is a local coordinate system on $M$, and $p= (p_1,\dots, p_n)$ are the corresponding momenta. 

Let $(S^n, g)$ be the standard sphere of dimension $n\ge 2$ in the Euclidean space $\mathbb{R}^{n+1}$ with standard induced metric.  We consider a (skew-symmetric) 2-form in $\mathbb{R}^{n+1}$ whose components are constants in Cartesian coordinates, and denote by $\omega=\omega_{ij}$
the restriction of this form onto the sphere.  We refer to $\omega$ as a constant magnetic form and  study the magnetic geodesic flow  on $S^n$ corresponding to $g$  and $\omega$. 

 \begin{theorem} 
 \label{thm:1} 
The magnetic geodesic flow on the sphere $(S^n,g)$ endowed with a constant magnetic form $\omega$ 
 is Liouville integrable by means of integrals linear and quadratic in momenta. 
 More specifically, there exist $n$ functions $F_1,\dots, F_n:T^*S^n\to \mathbb{R}$    such that  the following holds: 
 \begin{itemize}

 \item   $H$ is a linear combination of $F_1,\dots, F_n$ with constant coefficients. 

 \item  $F_1,\dots, F_n$  Poisson-commute with respect to the perturbed symplectic form $\Omega_{\mathrm{pert}}$.

 \item $F_1,\dots, F_n$   are functionally independent, i.e., their differentials are linearly independent almost everywhere.
 
\item The first  $\lfloor \frac{n}{2} \rfloor$ of the functions $F_1,\dots, F_n$ are quadratic in momenta with 
coefficients depending on the position,  the other $m=n-\lfloor  \frac{n}{2} \rfloor$ are linear.
 
 \end{itemize}
 \end{theorem}

 Theorem \ref{thm:1} proves  the conjecture recently proposed by Dragovic, Jovanovic and Gajic 
  \cite[Conjecture 5.1]{dragovic} and proved by them in dimension $n\le 5$ \cite{dragovic, dragovic1}.  Our approach is visually  different from that of \cite{dragovic, dragovic1} and is based on  some new ideas related to the  study of separation of variables and Killing tensors on constant curvature spaces \cite{BolsinovKonyaevMatveev-OrthogSep, Konrad1, Konrad2, Konrad3}, and in particular on separation of variables for the Neumann system.
  Independently, and almost simultaneously,  such integrability
has been shown by Dragovic, Jovanovic and Gajic in \cite{4}, where a Lax representation of the equations of motion was constructed. 

Our proof of Theorem \ref{thm:1} reduces the study of the magnetic geodesic flow on $(S^n, g,\omega)$ to the so-called degenerate Neumann system, which is known to be integrable by virtue of integrals linear and quadratic in momenta.  We show that the magnetic geodesic flow on $(S^n, g,\omega)$ can be equivalently formulated as the Hamiltonian system  on $T^*S^n$ with the canonical Poisson structure, whose Hamiltonian $H_{\mathrm{pert}} $ is the sum of the Hamiltonian of the degenerate Neumann system and a linear in momenta  integral which Poisson commutes with the integrals of the Neumann system.  Hence, the integrability automatically follows.

 \begin{remark}{\rm 
    If $\omega$ is a restriction  of the 2-form $\sum_{i,j=1}^{n+1} \alpha_{ij}\dd X^i\wedge \dd X^j$ with constant $\alpha_{ij}$
    onto the sphere $S^n=\Bigl\{ \sum_{k=1}^{n+1} (X^k)^2 = 1\Bigr\}\subset \mathbb R^{n+1}$, and  the matrix $\Bigl(\alpha_{ij}\Bigr)$  has multiple eigenvalues, then the magnetic geodesic flow is even superintegrable, in the sense that there exists additional integrals functionally independent of  $F_1,\dots, F_n$.  
 }\end{remark}

\subsection*{Acknowledgments.}  A.\,B. and A.\,K. were supported by the Ministry of Science and Higher Education of the Republic of Kazakhstan (grant No. AP23483476).   V.\,M. thanks the DFG (projects 455806247 and 529233771), and the ARC Discovery Programme DP210100951 for their support.

\section{Proof of Theorem  \ref{thm:1}}

\subsection{Reformulating the problem in terms of the canonical Poisson structure} \label{sec:2.1}
Let us
first  recall an equivalent  description  of the magnetic geodesic  flow on a manifold $M^n$ generated by a  metric $g$ and a closed 2-form $\omega$. 

Consider a 1-form $\sigma= \sigma_1\dd x^1+\cdots + \sigma_n\dd x^n$ such that $\dd \sigma=\omega. $ Local existence of $\sigma$ follows from the closedness  $\omega$. In  our setup, it  exists  globally, as the sphere is simply connected. \weg{We will give a formula for $\sigma$ later.}

Next, consider  the Hamiltonian system in the canonical Poisson structure on $T^*M$ generated by the Hamiltonian 
\begin{equation} \label{eq:1}
\begin{array}{rl}H_{\mathrm{pert}}:= &   \tfrac{1}{2}\sum_{i,j=1}^n g^{ij}(p_i- \sigma_i)(p_j-\sigma_j)\\=&     \tfrac{1}{2}\sum_{i,j=1}^n\left(g^{ij}p_ip_j - 2 \sigma_i g^{ij} p_j + g^{ij}\sigma_i\sigma_j\right).\end{array}
\end{equation}
\begin{fact}
    For any trajectory $(x^1(t),...,x^n(t), p_1(t),..., p_n(t))$ of this Hamiltonian  system,
the curve $$(x^1(t),...,x^n(t), p_1(t)+ \sigma_1(x(t)),..., p_n(t)+  \sigma_n(x(t)))$$
  is a trajectory of the magnetic geodesic flow, and vice versa.  
  Moreover, the transformation 
 \begin{equation} 
 \label{eq:transformation} 
 (x^1,...,x^n,p_1,...,p_n)\mapsto (x^1,...,x^n,p_1+ \sigma_1(x),...,p_n+\sigma_n(x)) 
 \end{equation}
maps $H_{\mathrm{pert}}$ to the unperturbed Hamiltonian \eqref{eq:Hunpert} and the canonical  symplectic form $\sum_{i=1}^n \dd p^i\wedge \dd x^i$ to the perturbed symplectic form \eqref{eq:pert}. 
In particular, it maps integrals of the Hamiltonian system generated by $H_{\mathrm{pert}}$, which Poisson commute with respect to the canonical symplectic structure, 
to integrals of the Hamiltonian system generated by \eqref{eq:Hunpert}, which Poisson commute with respect to the perturbed  symplectic structure \eqref{eq:pert}. \end{fact}

This fact is well known and is  easy to prove, as substituting $p_i+ \sigma_i$ instead of $p_i$ into  $\sum_{i=1}^n \dd p_i \wedge \dd x^i$ immediately gives the perturbed form 
\eqref{eq:pert}.

Thus,  instead of discussing integrability of the magnetic geodesic flow in its initial setup, i.e., in the sense of the perturbed symplectic form $\Omega_{\mathrm{pert}}$, in the proof of Theorem \ref{thm:1}, we may  and will discuss the integrability of the Hamiltonian system  generated by $H_{\mathrm{pert}}$ in the sense of the canonical form $\sum_{i=1}^n \dd p_i \wedge \dd x^i$. This viewpoint is more convenient, as it allows us to use known results
on separation of variables on the spaces of constant curvature  and  on integrability of  the (degenerate) Neumann system. Note that transformation \eqref{eq:transformation}  sends linear and quadratic in momenta functions  to (possibly, inhomogeneous) linear and quadratic in momenta functions.  

\subsection{Reduction to the degenerate  Neumann problem}
We consider the standard sphere $(S^n, g)$ of dimension $n\ge 2$
and a skew-symmetric 2-form on $\mathbb{R}^{n+1}$    whose entries are constant in the standard Cartesian coordinates 
 $X^1,\dots, X^{n+1} $.   
  As in the statement of Theorem \ref{thm:1}, we set  $m=n-\lfloor \frac{n}{2} \rfloor= \lfloor\frac{n+1}{2}\rfloor$. Without loss of generality, by \cite[end of \S 4 of Ch. XI]{Gantmacher},  we may assume that the constant form in 
$\mathbb{R}^{n+1}$ is  given by 
\begin{equation} \label{eq:2a}
\sum_{i=1 }^{m} \alpha_i \dd X^{2i-1}\wedge \dd X^{2i}     
\end{equation}
with nonnegative constants $ \alpha_i$. 

As a 1-form $\sigma$, in the ambient space $\mathbb{R}^{n+1}$, 
whose exterior derivative $\dd\sigma$ coincides with \eqref{eq:2a},   we choose  
\begin{equation} \label{eq:2b}
 \sum_{i=1 }^{m} \frac{1}{2}\left(X^{2i-1} \dd X^{2i} - X^{2i}\dd X^{2i-1}\right)\alpha_i . 
\end{equation}

Recall the exterior derivative of a form commutes with the restriction to a submanifold. The forms   \eqref{eq:2a} and \eqref{eq:2b}
are given in Cartesian coordinates $X^1,\dots, X^{n+1}$ in the ambient space $\mathbb{R}^{n+1}$; in order to obtain the forms  $\omega, \sigma$ 
on the sphere one should restrict them to the sphere. Let $\sigma = \sum_{i=1}^n \sigma_i \dd x^i$ denote this restriction in some local coordinates $x^1,\dots, x^n$ on $S^n$.  

Let us raise indices of $\sigma$, i.e., 
consider the vector field on the sphere    dual to the form \eqref{eq:2b}
with respect to the metric of the sphere. 
In order to obtain an expression for  it in the ambient coordinates,  
observe that  each form $X^{2i-1} \dd X^{2i} - X^{2i-1}\dd X^{2i}$  vanishes on the radial vector field  $X^1 \tfrac{\partial}{\partial X^1}+\cdots + X^{n+1}\tfrac{\partial }{\partial  X^{n+1}} $, which is orthogonal to the sphere. Then,  the {\it raising index} procedures for the ambient metric and for its restriction to the sphere coincide, and we obtain the vector field 
\begin{equation} 
\label{eq:2c} 
 \sum_{i=1 }^{m} \frac{1}{2}\left(X^{2i-1} \frac{\partial}{\partial  X^{2i}} -X^{2i}\frac{\partial }{\partial X^{2i-1}}\right)\alpha_i .
\end{equation}
 This vector field is tangent to the sphere, so its restriction to the sphere is well defined and coincides with $\sum_{j=1}^n\sigma^j\frac{\partial}{\partial x^j}$.  
Clearly, each term in the linear combination \eqref{eq:2c}, 
\begin{equation} 
\label{eq:2}
X^{2i-1} \frac{\partial}{\partial  X^{2i}} - X^{2i}\frac{\partial }{\partial X^{2i-1}}, 
\end{equation} 
is a Killing vector field, as it corresponds to the standard rotation in the plane with coordinates $X^{2i-1}, X^{2i}$.  Moreover,
these vector fields commute. 

Next, consider the terms $\sigma_i g^{ij} p_j$   and  $\tfrac{1}{2}g^{ij}\sigma_i\sigma_j$ from \eqref{eq:1}.  
We already know that  $\sigma^j= \sum_{i=1}^n \sigma_i  g^{ij}$ coincides with the  Killing vector field \eqref{eq:2c} so that $\sigma_i g^{ij} p_j$ is exactly the linear function on $T^*S^n$ corresponding to it. Next, the ``potential'' energy  $\tfrac{1}{2}\sum_{i,j=1}^n g^{ij}\sigma_i\sigma_j$  is just the scalar product of $\sigma^i$ with itself and, in the ambient coordinates,  
is the quadratic in $X^i$ function  
\begin{equation}\label{eq:potential} 
   \tfrac{1}{8} \sum_{i=1}^m \left((X^{2i-1})^2 + (X^{2i})^2\right)\alpha_i^2.  
\end{equation}

We see that the Hamiltonian \eqref{eq:1}, in our situation, is the  sum of the kinetic energy $K = \frac{1}{2}\sum_{ij} g^{ij}p_ip_j$ coming form the the standard metric of $S^n$, the potential energy \eqref{eq:potential} and the linear integral corresponding to the Killing vector field \eqref{eq:2c}. 

Now we note that the sum of the  kinetic energy $K$ and potential energy \eqref{eq:potential} gives the Hamiltonian  of the so-called degenerate Neumann system.  Recall that {\it Neumann system}  on $S^n$ is defined by the Hamiltonian $K + U$, where  $U$ is a quadratic potential of the form $\sum_{i=1}^{n+1} a_i(X^i)^2$ restricted to the sphere. A~Neumann system is {\it nondegenerate}, if  all the coefficients $a_i$ are different, and is {\it degenerate},  if  some of the them coincide. In our case, the Neumann system is degenerate, as 
the coefficients at $(X^{2i-1})^2$ and $(X^{2i})^2 $ are the same. Moreover, if certain constants $\alpha_i$ coincide, the ``level of degeneracy'' is higher, as more coefficients  coincide.  It is known that degenerate and nondegenerate Neumann systems are integrable in the class of quadratic in momenta integrals. For the nondegenerate system, the integrability was established e.g. in \cite{Babelon}.  For degenerate systems, see e.g. \cite{Dullin, Liu}.  In the next subsection, we will recall known results  about   nondegenerate  Neumann systems  (e.g. \cite{M83,M78}) and use them for describing the integrals of the degenerate Neumann problem which appears in our setting. The integrals should be chosen in such a way that they  Poisson commute with  the linear integral corresponding to the Killing vector field \eqref{eq:2c}. 

\subsection{ Uhlenbeck integrals for the Neumann system, and integrability for certain degenerate Neumann systems.  } \label{sec:2.3}

We consider the Neumann problem of a point moving on the sphere 
$$
S^{n} =\bigl\{ (X^1,\dots, X^{n+1})\in \mathbb{R}^{n+1} \mid (X^1)^2 +\dots + (X^{n+1})^2  = 1\bigr\}
$$
under a quadratic potential 
$$
U_A =  a_1 (X^1)^2 +\dots + a_{n+1} (X^{n+1})^2 .
$$
We think of it as a Hamiltonian system on $T^*S^n$.   

We use the following notation $M_{ij}=X^i \tfrac{\partial }{\partial X^j} - X^j\tfrac{\partial}{\partial X^i}$ for the standard basis in the space of Killing vector fields or, equivalently, in the isometry Lie algebra $so(n+1)$.  Notice that we may think of $M_{ij}$ as a linear function on the cotangent bundle $T^*S^{n}$, so that the expression $M_{ij}^2$ below is understood as an {\it elementary} quadratic function on $T^*S^n$.  In this notation, the Hamiltonian of the Neumann problem takes the form
\begin{equation} \label{eq:neumann}
H = K + U_A, \quad \mbox{where} \ K=\frac{1}{2}\sum_{i<j} M_{ij}^2.
\end{equation}

The integrability in the generic case, when   all $a_i$ are different, is established by the following  well known result. 

\begin{fact}[e.g., \cite{M83, M78}]
\label{fact:2}
 Let $a_i\ne a_j$ for $i\ne j$. Then Poisson commuting integrals of the Neumann problem can be taken in the form
\begin{equation}
\label{eq:9}
F_B = K_B + U_B,  \quad \mbox{with $B=(b_1,\dots,b_{n+1})\in\mathbb{R}^{n+1}$},
\end{equation}
where 
\begin{equation}
\label{eq:10}
K_B = \frac{1}{2} \sum_{i<j} \frac{b_i-b_j}{a_i-a_j}M_{ij}^2,   \quad U_B = \sum b_i x_i^2. 
\end{equation}

The integrals $F_{B_1},\dots,F_{B_n}$ are functionally independent if and only if the vectors $B_1,\dots, B_n$ are linearly independent and  $(1,1,\dots,1)$ does not belong to $\operatorname{Span}(B_1,\dots,B_n)$.  In particular, these integrals guarantee Liouville integrability of the nondegenerate  Neumann problem. 
\end{fact}

In \cite{M83, M78},  the integrals \eqref{eq:9}, written in a slightly different but equivalent form, were attributed to K.\,Uhlenbeck.

\begin{remark}\label{rem1}{\rm It follows from the above formulas that the collection of functions $\{ F_B, \ B\in \mathbb{R}^{n+1}\}$ is a vector space of dimension $n+1$.  Indeed,  $F_{\lambda_1 B_1+\lambda_2B_2}=\lambda_1 F_{B_1}+\lambda_2 F_{B_2}$, and moreover $F_B=0$ if and only if $B=0$.  However, $F_{(1,\dots,1)}=\sum x_i^2 =1$ is a constant function on  $T^*S^n$, which should be treated as  a  trivial/ignorable integral.   The functions $F_{B_1},\dots, F_{B_n}$ from the last statement of Fact  \ref{fact:2} can be naturally understood as a {\it basis  of $\{ F_B, \ B\in \mathbb{R}^{n+1}\}$ modulo constants}.  Fact  \ref{fact:2} basically says that the functions from such a basis are not only linearly, but also functionally independent on $T^*S^n$. 

Note also that for $B=A$ the integral  $F_B$ is the Hamiltonian of the Neumann system.

For our purposes, we will also need to deal with the homogeneous quadratic parts $K_B$ of functions $F_B$.   
Notice that $K_B=0$ if and only if $B=(\lambda,\dots,\lambda)$ so that $\dim \{K_B, \ B\in\mathbb{R}^{n+1}\}=n$ and every basis $K_{B_1},\dots,K_{B_n}$ of 
$\{K_B, \ B\in\mathbb{R}^{n+1}\}$ provides $n$ Poisson commuting independent integrals of the geodesic flow on $S^n$. Moreover, 
at almost every point  $x\in T^*S^n$, there exists a basis in $T^*_xS^n$, such that in this basis the matrices of  all $K_B$ are diagonal. 
}\end{remark}

Next, consider the case when $A$ is singular in the sense that some of $a_i$ coincide:
\begin{equation}
\label{eq:Asing1}
a_1=\dots=a_{k_1} < a_{k_1+1} = \dots = a_{k_1+k_2} < \dots < a_{k_1+\dots+k_{s-1}+1} = \dots = a_{k_1+\dots+k_s}
\end{equation}
In other words,  the collection of indices $\{1,2,\dots, n+1\}$ is partitioned into $s$ subsets $I_1,\dots, I_s$.  The $r$-th subset consists of $k_r$ indices that correspond to equal $a_i$'s, more specifically, 
{\small
$$
I_r = \Bigl\{k_1+\dots+k_{m-1} +1, \dots, k_1 + \dots +k_r\Bigr\} \quad \mbox{and} \quad  k_1+k_2+\dots+k_s = n+1.
$$} 
 
For our further purposes,  consider \begin{equation} \label{eq:Gm}
\mathcal G_r = \operatorname{Span}\Bigl( M_{lm},  \ l,m\in I_r\Bigr).\end{equation}   Obviously $\mathcal G_r$ is a subalgebra of the algebra of Killing vector fields, which is isomorphic to $so(k_r)$.  

\weg{
The {\it commuting} first integrals we are going to construct will be quadratic forms in $M_{ij}$. For this reason, together with $\mathcal G_r$ we consider its {\it quadratic} extension
$$
{\mathcal G}^2_r = \{ \textrm{all quadratic, not necessarily homogeneous, polynomials in  $M_{ij}\in{\mathcal G}_r$}
\} 
$$
}

For a given $A$, we introduce the collection of (non-homogeneous) quadratic functions  $\mathcal F_A$ of the form
\begin{equation}
\label{eq:F_Bsing}
F_B = K_B + U_B
\end{equation}
with
\begin{equation}
\label{eq:K_Bsing}
K_B = \frac{1}{2} \sum_{i<j, \, a_i\ne a_j} \frac{b_i-b_j}{a_i-a_j}M_{ij}^2,   \quad U_B = \sum b_i x_i^2, 
\end{equation}
\begin{equation}
\label{eq:Bsing}
B=(b_1,\dots,b_{n+1})\in\mathbb{R}^{n+1}, \quad
b_1=\dots=b_{k_1}, \  b_{k_1+1} = \dots = b_{k_1+k_2},  \  \dots \  
\end{equation}
Notice that the components $b_i, b_j$ of $B$ are equal if $a_i=a_j$; but $b_i$ may be equal to $b_j$ even if $a_i\ne a_j$. 

As compared to formulas for $K_B$ in Fact
\ref{fact:2},  we simply remove all the terms which contain division by zero. The collection of quadratic functions $K_B$ defined by \eqref{eq:K_Bsing} (i.e., obtained from $\mathcal F_A$ by removing potentials)  will be denoted by $\mathcal K_A$.

The matrix $B$ in \eqref{eq:Bsing} depends on $s$ free parameters. The same argument as in Remark \ref{rem1} shows that $\dim \mathcal F_A=s$, $\dim \mathcal K_A=s-1$ and $F_{(1,\dots,1)}=1$.  Moreover,  if $B_1, \dots, B_{s-1}$ are vectors as in \eqref{eq:Bsing} which are linearly independent modulo $B=(1,\dots,1)$, then the functions $F_{B_1},\dots, F_{B_{s-1}}$ form a basis of  $\mathcal F_A$ modulo constants.   Similarly, their quadratic parts $K_{B_1},\dots, K_{B_{s-1}}$ form a basis of $\mathcal K_A$. Note also that the function $F_A$ is the Hamiltonian of the (degenerate) Neumann system.

Let us emphasise that $\mathcal F_A$  (as well as $\mathcal K_A$) is a well defined collection of functions for any $A\in\mathbb{R}^{n+1}$ satisfying \eqref{eq:Asing1}.  In particular, if all the components of $A$ are different, we obtain exactly the collection of functions from Fact \ref{fact:2}.  Also notice that the condition that the components of $A$ are arranged in ascending order is made only for convenience.  The construction can be naturally reformulated for an arbitrary $A$.

From Fact \ref{fact:2} we can easily derive the following statement. 

\begin{corollary} \label{cor:1}
For a fixed partition $I_1,\dots,I_s$, consider $B_1$ and $B_2$ satisfying \eqref{eq:Bsing}.
Then $F_{B_1}$ and $F_{B_2}$  Poisson commute. 

Moreover,  any element of $\mathcal G_r$,  $r=1,\dots, s$,  Poisson commutes with $F_{B_1}$ and $F_{B_2}$. Furthermore,  any  element of  $\mathcal G_{r_1}$ Poisson commutes with any element of $\mathcal G_{r_2}$ for  $r_1\ne r_2$, $r_1,r_2\in \{1,\dots, s\}$.  
\end{corollary}

Of course, if $k_r\ge 3,$ the elements of $\mathcal G_{r}$ do not commute, as the algebra $so(k_r)$ is not commutative.

\begin{proof}
The second statement of Corollary is obvious, as both the kinetic and potential parts of the  function
$F_B$ are preserved by the flows of the Killing vector fields $M_{ij}\in \mathcal G_r$. The third statement is also trivial, as the components from different $\mathcal G_r$'s depend on different groups of coordinates. 

In order to prove the first statement, we use the `passage to limit' procedure. We consider a converging  sequence $A(1),A(2),\dots, A(\ell),\dots \ \stackrel{\ell\to \infty}{\longrightarrow} A$, such that $A(\ell)$ is nonsingular, in the sense that  all of its entries are different. 
    
    Next, consider the integrals $F_{B_1}(\ell)  $ and $F_{B_2}(\ell)$ constructed by $A(\ell)$  and by $B_1$ and $B_2$. 
    We assume that $B_1$ and $B_2$ satisfy \eqref{eq:Bsing}. The functions $F_{B_1}(\ell)  $ and $F_{B_2}(\ell)$ Poisson commute, for every $\ell$, and converge to the integrals $F_{B_1} $ and $F_{B_2}$ as $\ell \to \infty$. Passing to the limit, we obtain the desired statement. 
 \end{proof}

The special case when each $I_r$ has at most two elements  is especially  important to our initial problem.

\begin{corollary} \label{cor:2}
 Assume that the entries $a_i$ of  $A$  satisfy 
\begin{equation}
    \label{eq:12}
    a_1=a_2<a_3=a_4<a_5=a_6<\cdots \  . 
\end{equation}
(if $n$ is even, the sequence of equalities and inequalities  \eqref{eq:12}
 ends  as follows: \   $\cdots  =a_{n}<a_{n+1}$. If $n$ is odd, it ends with  $\cdots = a_{n-1}<a_n=a_{n+1}$).

 Then the collection consisting of the quadratic functions $F_B$ defined by \eqref{eq:F_Bsing}--\eqref{eq:K_Bsing}  and linear functions $M_{2i-1,2i}$,  $1\le i \le \lfloor\frac{n+1}{2}\rfloor$,  is 
  Poisson commutative. 
 \end{corollary}

Under the assumptions of Corollary \ref{cor:2},  
the functions   $M_{2i-1, 2i}$   are clearly  functionally independent of the functions $F_B$. The number of the functions $M_{2i-1, 2i}$  is  $m=\lfloor\frac{n+1}{2}\rfloor$, and the number of  functionally independent functions  $F_B$ is  $n-m=\lfloor\frac{n}{2}\rfloor$, so these integrals insure the Liouville integrability of the (degenerate) Neumann problem with the potential $U_A$ (recall that the Hamiltonian of the Neumann system is $F_B$ with $B=A$). Note also that $F_B$'s are simultaneously diagonalisable in a certain basis at almost every point of the sphere.

\begin{remark}{\rm
Corollary \ref{cor:2}  proves Theorem \ref{thm:1} under the additional assumption that the magnetic form $\omega$ is the restriction of the form \eqref{eq:2a} with $\alpha_i\ne \alpha_j$ for $i\ne j$. Indeed, the perturbed Hamiltonian $H_{\mathrm{pert}}$ is obtained from the Hamiltonian $H=K+U_A$ of the Neumann problem by adding the linear function corresponding to the vector field \eqref{eq:2c}, that is,
$$
H_{\mathrm{pert}}=H+\tfrac{1}{2} \left(  \alpha_1 M_{12} + \alpha_2 M_{34} + \alpha_3 M_{56} + \dots  \right)
$$
Thus, the integrals from Corollary \ref{cor:2} Poisson commute with  $H_{\mathrm{pert}}$ and, therefore, guarantee Liouville integrability by means of quadratic and linear integrals as stated in Theorem \ref{thm:1}. 
One can also show that the above integrals naturally lead to separation of variables in the sense of St\"ackel.  
}\end{remark}

\subsection{The   existence  of  integrals  commuting with $M_{2i-1, 2i}$ in the general case.}

We now allow some of the constants  $\alpha_i$ to be equal. Our goal is to show the existence of sufficiently many   quadratic in momenta integrals, 
commuting with the integrals $M_{2i-1, 2i}$ coming from the Killing vector fields \eqref{eq:2}. In \S \ref{sec:2.3}, we did this under the assumption that   all $\alpha_i$'s are different. 

The general case will be done by the passing to limit procedure: we consider $m$ 
sequences  
\begin{equation} 
\label{eq:sec} 
k\mapsto \alpha_1(k), \ k\mapsto \alpha_2(k), \dots, \ k\mapsto \alpha_{m}(k),
\end{equation} 
such that for any $k$ we have $\alpha_i(k)\ne \alpha_j(k)$ for $i\ne j$,  $\alpha_i(k)$ are all nonnegative and such that $\lim_{k\to \infty} \alpha_i(k)= \alpha_i$. 
By Corollary \ref{cor:2}, for 
each  $k$,  there exists an $n$-dimensional space generated by $n$ functionally independent  quadratic in momenta  integrals\footnote{Strictly speaking, Corollary \ref{cor:2} provides independent integrals some of which are linear. To get a collection of quadratic integrals, we can just square them.} 
$$
\textrm{Span}\Bigl\{F_1(k)=K_1(k) + U_1(k), \dots, F_n(k)=K_n(k)+ U_n(k)\Bigr\}
$$ 
such that any element of this space is 
invariant with respect to the Killing vector fields \eqref{eq:2c}.    Without loss of generality, we assume that $F_1(k)$ is the Hamiltonian of the Neumann system corresponding to $\alpha_1(k),\dots, \alpha_m(k)$. 

In order to define the limit of such spaces of integrals, we will first define the limit of  the space of their ``kinetic'' parts $K_i$. We employ the approach developed and used by K. Sch\"obel at al, see e.g.  \cite{Konrad1,Konrad2, Konrad3}. By \cite{smirnov},  to each homogeneous quadratic in momenta integral of  the geodesic flow on $S^n$, one can canonically, by a real-analytic  formula, assign a tensor $R_{IJKL}$ on  $\mathbb{R}^{n+1}$
satisfying  the symmetries of the curvature tensor, whose entries are constants in the ambient coordinates 
$X^1,\dots, X^{n+1}$.    We denote  the space of such $(0,4)$ tensors  
by $\mathbf{K}$.  

The corresponding mapping $\phi$ from the space of homogeneous quadratic integrals to $\mathbf{K}$ is a linear isomorphism. We emphasise that the tensor  $\phi(Q)= R_{IJKL}$ ``knows everything'' about the homogeneous quadratic integral $Q$. In particular, the entries of $Q$ and their derivatives can be reconstructed by $R_{IJKL}$ by an algebraic procedure. 

In our situation,  the sequences \eqref{eq:sec} gives us a sequence of  
$n$-dimensional vector subspaces in the space of  quadratic integrals.  Combining it with $\phi$, we obtain  a sequence of   $n$-dimensional vector subspaces of $\mathbf{K}$.  Since the space of $n$-dimensional  vector subspaces of $\mathbf{K}$ is evidently compact, 
  the sequence  has a convergent subsequence.  Without loss of generality, we think  that  the initial sequence converges.  The  limit is then an $n$-dimensional subspace of  $\mathbf{K}$. As  $\phi$ is a bijection, we obtain 
  an $n$-dimensional space of quadratic in momenta functions which are integrals for the geodesic flow on $S^n$.

  Next, observe that the Poisson commutativity for quadratic integrals is an algebraic  condition on  the entries of the integrals and their first derivatives. Hence, it is an algebraic condition on the entries of the corresponding elements of $\mathbf{K}$. 
  As this condition was fulfilled for all  elements of the sequence, it is fulfilled for the limit as well. We therefore obtain an $n$-dimensional linear family of  Poisson  {\it commuting} integrals of the geodesic 
 flow on $S^n$. We denote a basis in this family by $K_1,\dots, K_n$, thinking of $K_1$ as the kinetic energy of the standard metric on $S^n$.   
 
 The integrals corresponding to $\alpha_1(k), \dots,  \alpha_m(k)$ were, 
 by construction,  invariant with respect to the Killing vector fields \eqref{eq:2}. Then, the quadratic functions  $K_i$  are also invariant with respect to these Killing vector fields, and therefore with respect to the Killing vector field \eqref{eq:2c}. 

Note also that  $K_1(k),\dots, K_n(k)$ are simultaneously diagonalisable, at almost every point $x\in S^n$, in a certain frame in $T^*S^n$. Passing to the limit, we obtain that the integrals  $K_1,\dots, K_n$ are also simultaneously diagonalisable. Then,
linear independence implies functional independence of $K_1,\dots, K_n$. 

 Let us now add potential energies to the construction. 
First observe that for two Poisson commuting homogeneous quadratic 
functions $F_1= \sum_{i,j=1}^n K^{ij}p_ip_j $ and $F_2= \sum_{i,j=1}^n L^{ij}p_ip_j $, the condition that   $F_1+ U$ and $F_2 + V$ Poisson commute is equivalent to the relation 
\begin{equation} 
\label{eq:14}
\sum_{s=1}^n K^{si}\tfrac{\partial V}{\partial x^s} = \sum_{s=1}^n L^{si}\tfrac{\partial U}{\partial x^s} . 
\end{equation}
 If the  kinetic part of the integral corresponds to the metric, i.e., $K^{ij}=  g^{ij}$,  then
 the necessary and sufficient condition for local existence of a function $V$, satisfying  \eqref{eq:14}  for a given $U$, is the so-called {\it Benenti condition} 
   \begin{equation}
   \label{eq:15}  
   \dd\left(\sum_{s,i=1}^nL^s_i \frac{\partial U}{\partial x^s}  \dd x^i\right)=0, 
   \end{equation} where we used $g$ for index manipulations.
 Moreover, such a function $V$, if exists,  is unique  up to adding a constant and satisfies the equation 
 \begin{equation}
 \label{eq:16} 
 \dd V =  \sum_{s,i=1}^nL^s_i \frac{\partial U}{\partial x^s} \dd x^i.
 \end{equation}
Note that the sphere is simply connected, so if \eqref{eq:15} is fulfilled, then there exists a global solution of \eqref{eq:16}.

In our setting, the sequence of potential energies   $U_{1}(k)$ of the  Neumann systems   corresponding to the constants $\alpha_1(k), \dots,  \alpha_m(k)$, 
evidently converges to the  potential energy of the Neumann system corresponding to  
$\alpha_1, \dots,  \alpha_m$. As in the nondegenerate case,  each function $U_1(k)$ satisfies  the Benenti condition \eqref{eq:15} with respect to each $K_i(k)$. 
Passing to the limit, we obtain that 
the potential energy $U_1$ of the Neumann system  corresponding to 
$\alpha_1, \dots,  \alpha_m$ satisfies the Benenti condition \eqref{eq:15} with respect to the quadratic parts  $K_1, \dots, K_n$. 
As the sphere is simply-connected,  there exist functions $U_i$ such that  $K_i+ U_i$ Poisson commute  with the  the Hamiltonian of our Neumann system. 
Note that since \eqref{eq:16} is invariant with respect to the flows of the vector fields \eqref{eq:2}, the  functions $U_i$ are invariant with respect to them also.

In order to show that  $F_i=K_i+ U_i$   Poisson commute pairwise, 
we use the fact that  for any $k$ the functions constructed by the formula \eqref{eq:16} with $L= K_i(k)$ and $U=U_1(k)$  are, up to constants, the potential parts of the integrals $F_i(k)=K_i(k)+U_i(k)$ of the  Neumann system corresponding to $\alpha_1(k),\dots, \alpha_m(k)$.  Then these functions satisfy  \eqref{eq:14}. Passing to the limit, we obtain that the functions $U_k$ also satisfy relation \eqref{eq:14} and therefore the corresponding integrals $F_i$ Poisson commute.

Clearly, the functions $F_i$ are functionally independent, as their quadratic parts are functionally independent. We have shown above that they are invariant with respect to the flows of the Killing vector fields \eqref{eq:2} and therefore commute with the corresponding  integrals $M_{2i-1,2i}$ linear in momenta. They also commute with the linear integral corresponding to the vector field \eqref{eq:2c}, and therefore one can replace, keeping the integrability, 
the last $m$ integrals by the linear integrals $M_{2i-1,2i}$. 

Thus, we have shown that the existence of $n$ Poisson commuting functionally independent 
functions $F_1,\dots, F_n$ such that the first $n-m$ are quadratic in momenta, the last $m$ are linear in momenta,  and the Hamiltonian $H_{\mathrm{pert}}$ given by \eqref{eq:pert} is their linear combination. Theorem \ref{thm:1} is proved. 

{
Notice that the above `passage to limit'  construction is quite general and can be applied to various integrable systems depending on parameters when one needs to study their degenerations.  Alternatively, in our case this passage to limit can be made very explicit. Indeed, consider $A=(a_1,\dots,a_{n+1})$ and $B=(b_1,\dots,b_{n+1})$  as in \eqref{eq:Asing1}, \eqref{eq:Bsing},  and choose the deformations  $A(t)\to A$, $B(t)\to B$ as follows:
$$
a_i(t) = a_i + t\lambda_i \quad \mbox{and} \quad b_i(t) = b_i + t\mu_i, 
$$
To put everything into the context of magnetic flows, we assume in addition that
$a_{2i-1}(t)\equiv a_{2i}(t)$, $b_{2i-1}(t)\equiv b_{2i}(t)$.
Then the integrals $F_{B(t)}$ from \eqref{eq:F_Bsing} take the form
$$
F_{B(t)}=\tfrac{1}{2}\sum_{a_i(t)\ne a_j(t)} \tfrac{b_i-b_j + t(\mu_i-\mu_j)}{a_i-a_j + t(\lambda_i-\lambda_j)} M_{ij}^2 +
\sum (b_i + t\mu_i) (X^i)^2
$$
and the passage to limit as $t\to 0$ can be easily performed for each term separately, as no division by zero appears if $\lambda_i$'s are appropriately chosen. Since the parameters $b_i$ and $\mu_i$ are  free and independent of each other, we obtain a collection of commuting quadratic integrals of two types:
$$
F_B = \tfrac{1}{2} \sum_{a_i \ne a_j} \tfrac{b_i-b_j}{a_i-a_j} \, M_{ij}^2  + \sum b_i (X^i)^2\quad \mbox{as in Corollary \ref{cor:1},}
$$
and
$$
F_{I_r , \mu}=\tfrac{1}{2} \sum_{l,m\in I_r, \lambda_l\ne\lambda_m} \tfrac{\mu_l-\mu_m}{\lambda_l-\lambda_m} \, M_{lm}^2 
$$
The latter is a quadratic form in the generators $M_{lm}$ of the subspace $\mathcal G_r$ defined in \eqref{eq:Gm}.  

The number of independent integrals of the form $F_{B(t)}$ for $t\ne 0$ equals $\lfloor\frac{n}{2}\rfloor$  (see Corollary \ref{cor:2}). One can check that the above collection  still contains the same number of independent quadratic integrals. Moreover, before and after taking the limit,  all these functions commute with $m=\lfloor\frac{n+1}{2}\rfloor$ linear functions $M_{2i-1,2i}$ and, therefore, with the linear function associated with the vector field \eqref{eq:2c}, as required. 
}

\printbibliography

\end{document}